\documentclass[a4paper,12pt]{amsart}
\usepackage{amssymb}
\usepackage{ifthen}
\nonstopmode \numberwithin{equation}{section}
\newtheorem{thm}{Theorem}
\newtheorem{cor}{Corollary}
\newtheorem{lem}{Lemma}

\newtheorem{conj}{Conjecture}

\theoremstyle{definition}
\newtheorem{defn}{Definition}[section]

\newtheorem{prob}[equation]{Problem}

\newenvironment{rem}{%
\bigskip
\noindent \textsl{{\sl Remark. }}}{\bigskip}
\newenvironment{rems}{%
\bigskip
\noindent \textsl{{\sl Remarks. }}}{\bigskip}

\newcounter {own}
\def\theown {\thesection       .\arabic{own}}

\newenvironment{pf}[1][]{%
 \vskip 3mm
 \noindent
 \ifthenelse{\equal{#1}{}}%
  {{\slshape Proof. }}%
  {{\slshape #1.} }%
 }%
{\qed\bigskip}

\newcounter{alphabet}
\newcounter{tmp}
\newenvironment{Thm}[1][]{\refstepcounter{alphabet}%
\bigskip%
\noindent%
{\bf Theorem \Alph{alphabet}}%
\ifthenelse{\equal{#1}{}}{}{ (#1)}%
{\bf .} \itshape}{\vskip 8pt}

\newcommand{\ID}{{\mathbb D}}

\newcommand{\IC}{{\mathbb C}}

\newcommand{\sphere}{{\widehat{\mathbb C}}}

\def\be{\begin{equation}}
\def\ee{\end{equation}}

\newcommand{\bee}{\begin{enumerate}}
\newcommand{\eee}{\end{enumerate}}

\newcommand{\blem}{\begin{lem}}
\newcommand{\elem}{\end{lem}}
\newcommand{\bthm}{\begin{thm}}
\newcommand{\ethm}{\end{thm}}
\newcommand{\bcor}{\begin{cor}}
\newcommand{\ecor}{\end{cor}}
\newcommand{\beg}{\begin{examp}}
\newcommand{\eeg}{\end{examp}}
\newcommand{\begs}{\begin{examples}}
\newcommand{\eegs}{\end{examples}}
\newcommand{\bdefe}{\begin{defn}}
\newcommand{\edefe}{\end{defn}}
\newcommand{\bprob}{\begin{prob}}
\newcommand{\eprob}{\end{prob}}
\newcommand{\bei}{\begin{itemize}}
\newcommand{\eei}{\end{itemize}}

\newcommand{\bcon}{\begin{conj}}
\newcommand{\econ}{\end{conj}}
\newcommand{\bcons}{\begin{conjs}}
\newcommand{\econs}{\end{conjs}}
\newcommand{\bprop}{\begin{propo}}
\newcommand{\eprop}{\end{propo}}
\newcommand{\br}{\begin{rem}}
\newcommand{\er}{\end{rem}}
\newcommand{\brs}{\begin{rems}}
\newcommand{\ers}{\end{rems}}
\newcommand{\bo}{\begin{obser}}
\newcommand{\eo}{\end{obser}}
\newcommand{\bos}{\begin{obsers}}
\newcommand{\eos}{\end{obsers}}
\newcommand{\bpf}{\begin{pf}}
\newcommand{\epf}{\end{pf}}
\newcommand{\ba}{\begin{array}}
\newcommand{\ea}{\end{array}}
\newcommand{\beq}{\begin{eqnarray}}
\newcommand{\beqq}{\begin{eqnarray*}}
\newcommand{\eeq}{\end{eqnarray}}
\newcommand{\eeqq}{\end{eqnarray*}}

\newcommand{\ra}{\rightarrow}

\newcommand{\ds}{\displaystyle}

\newcounter{minutes}
\divide\time by 60
\newcounter{hours}
\multiply\time by 60 \addtocounter{minutes}{-\time}

\begin{document}

\title{Criteria for univalence, Integral means and Dirichlet integral for Meromorphic functions}
\begin{center}
{\tiny \texttt{FILE:~\jobname .tex,
        printed: \number\year-\number\month-\number\day,
        \thehours.\ifnum\theminutes<10{0}\fi\theminutes}
}
\end{center}

\author{Bappaditya Bhowmik${}^{~\mathbf{*}}$}
\address{Bappaditya Bhowmik, Department of Mathematics,
Indian Institute of Technology Kharagpur, Kharagpur - 721302, India.}
\email{bappaditya@maths.iitkgp.ernet.in}
\author{Firdoshi Parveen}
\address{Firdoshi Parveen, Department of Mathematics,
Indian Institute of Technology Kharagpur, Kharagpur - 721302, India.}
\email{frd.par@maths.iitkgp.ernet.in}

\subjclass[2010]{30C45, 30C70}
\keywords{ Meromorphic function, Concave function, Starlike function, Dirichlet finite integral, Integral mean }
\begin{abstract}
Let $\mathcal{A}(p)$ be the class consisting of functions $f$ that are holomorphic
in $\ID\setminus \{p\}$, $p\in (0,1)$ possessing a simple pole at the point $z=p$ with nonzero residue and normalized by the condition $f(0)=0=f'(0)-1$.
In this article, we first prove a sufficient condition for univalency for functions in $\mathcal{A}(p)$.
Thereafter, we consider  the class denoted by $\Sigma(p)$ that consists of functions $f \in \mathcal{A}(p)$ that are  univalent in
$\ID$. We obtain the exact value for $\ds\max_ {f\in \Sigma(p)}\Delta(r,z/f)$, where the Dirichlet integral $\Delta(r,z/f)$ is given by
$$
\Delta(r,z/f)=\ds\iint_{|z|<r} |\left(z/f(z)\right)'|^2 \,dx\, dy, \quad(z=x+iy),~0<r\leq 1.
$$
We also obtain a sharp estimate for $\Delta(r,z/f)$ whenever $f$ belongs to certain subclasses of $\Sigma(p)$.
Furthermore, we obtain sharp estimates of the integral means for the aforementioned classes of functions.

\end{abstract}
\thanks{}
\maketitle
\pagestyle{myheadings}
\markboth{B. Bhowmik and F. Parveen}{Criteria for univalence, Integral means and Dirichlet integral for Meromorphic functions}

\section{Introduction}
We use the following notations throughout the discussion of this article. Let $\ID:=\{z\in \mathbb{C}:|z|<1\}$ be the open unit disc where
$\IC$ is the whole complex plane. Let $\sphere$ denote the set $\IC \cup\{\infty\}$. We now recall the following basic classes of  functions
which are the main objects of study of many function theorists for several years now. Let $\mathcal{H}$ be the family of analytic functions in $\ID$ and
$\mathcal{A}$ be the subfamily of  $\mathcal{H}$ consisting of functions $g$ that satisfy the normalization $g(0)=0=g'(0)-1$. We consider the class
$\mathcal{S}:=\{g\in \mathcal{A}: g \mbox{ is univalent in }  \ID\}$.
Clearly $\mathcal{S}\subsetneq \mathcal{H}$. Let $\mathcal{C}$ and
$\mathcal{S}^{*}$ be the subclasses of $\mathcal{S}$ which are convex ($f(\ID)$ is a convex set) and starlike ($f(\ID)$ is a starlike set with respect to
the origin) respectively. We also consider the class $\Sigma$ of meromorphic
univalent functions in $\ID^{*}:=\{z\in \sphere :|z|>1\}$ having a simple pole at infinity with residue $1$.
We now discuss about the motivation and background of the problems that we consider in this article.
Let  $g\in \mathcal{H}$. We denote the area of the image of the disk $\ID_r:=\{z: |z|<r\}$
under $g$ by $\Delta(r,g)$, where $0<r\leq1$ and
$$
\Delta(r,g):=\iint_{\ID_r} |g'(z)|^2 \,dx\, dy, \quad(z=x+iy).
$$
The above integral $\Delta(r,g)$ is popularly known as {\it Dirichlet integral}. Each function $g\in \mathcal{H}$ has the Taylor expansion
$g(z)=\sum_{n=0}^{\infty}a_nz^n $ in $\ID$ and consequently, we have $g'(z)=\sum_{n=1}^{\infty}na_nz^{n-1}$.  It is now a simple exercise to compute
\beq \label{fp3eq1}
\Delta(r,g)=\pi \sum_{n=1}^{\infty}n|a_n|^2r^{2n}.
\eeq
Moreover, if $g\in \mathcal{S}$, we have $a_0=0$, $a_1=1$ and
$$
\Sigma \ni T(g)(z)= g(1/z)^{-1}=z-a_2+\sum_{n=1}^{\infty}c_nz^{-n},\quad z\in \ID^{*}.
$$
Now an application of Gronwall's area theorem applied to the above function $T(g)$ will yield $\sum_{n=1}^{\infty}n|c_n|^2\leq 1$.
For $g\in \mathcal{S}$, we have the following expansion for $z/g$:
$$
z/g(z)= 1-a_2z+\sum_{n=1}^{\infty}c_nz^{n+1},\quad z\in \ID.
$$
Now considering the above form of $z/g$, an application of the Gronwall's area inequality ($\sum_{n=1}^{\infty}n|c_n|^2\leq 1$) along with
the fact that $|a_2|\leq 2$ , S. Yamashita (compare \cite[Theorem 1]{Yama}) obtained:

\begin{Thm}\label{Yamashita} For $g\in \mathcal{S}$, we have
$$
\ds\max_ {g\in \mathcal{S}}\Delta(r,z/g)= 2\pi r^2 (r^2+2), \quad 0< r \leq 1.
$$
For each $r\in (0, 1]$, the maximum is attained only by the rotation of the Koebe function $K_\theta(z)=z/(1-e^{i\theta}z)^2$, $\theta \in (0, 2\pi]$.
\end{Thm}

In the same article (compare \cite[p. 438]{Yama}), Yamashita conjectured that
$$
\ds\max_ {g\in \mathcal{C}}\Delta(r,z/g)= \pi r^2, \quad 0< r \leq 1,
$$
where the maximum is attained only by the rotations of the function $g(z) = z/(1-e^{i\theta}z)$, $\theta \in (0, 2\pi]$.
This conjecture has recently been settled by M. Obradovic et.al. in \cite{OPW-1}. In a recent article (see \cite{PY}), Ponnusamy and Abu Muhanna
have obtained sharp estimates for the generalized Yamashita functional i.e. $\Delta(r,\phi(z)/f(z))$ for the class
of concave univalent functions with opening angle $\pi\alpha$, $\alpha\in (1,2]$ at infinity,  where $\phi$ is a Schwarz function.

In this article, we would like to consider meromorphic univalent functions with pole at $z=p\in (0,1)$.
Let $\mathcal{A}(p)$ be the class consisting of functions $f$ that are holomorphic
in $\ID\setminus \{p\}$, possessing a simple pole at the point $z=p$ with nonzero residue and normalized by the condition $f(0)=0=f'(0)-1$.
Let $\Sigma(p):=\{f\in\mathcal{A}(p): f \mbox{ is one to one in} ~~\ID \}$. We organize the paper as follows.
In the next Section, i.e. in Section 2, before we present our main results,  we first establish a sufficient condition for univalence for functions
in $\mathcal{A}(p)$ and we feel that it will be useful to present the absolute estimates
for the Dirichlet integrals $\Delta(r,f)$ and $\Delta(r,f/z)$
for $f\in\Sigma(p)$ and $0<r<p$. We also verify that these
results coincide with those of Yamashita in \cite{Yama} for the analytic case as we take the limit $p\ra 1^{-}$.

Another interesting subclass of $\Sigma(p)$
has recently been introduced by the authors of the present article in \cite{BF-1}. This class is denoted by $\mathcal{U}_{p}(\lambda)$
and  consists of all functions $f \in \mathcal{A}(p)$ such
that $\left|U_f(z)\right|< \lambda\mu$ for some $0<\lambda \leq1$ where
$$
U_f(z):=\left(z/f(z)\right)^2f'(z)-1 \quad \mbox{and} ~\mu:=\left((1-p)/(1+p)\right)^2.
$$
It has been shown in \cite{BF-1} that $\mathcal{U}_{p}(\lambda) \subsetneq \Sigma(p)$ and the interested reader may look at this article
for many other results on this newly defined class of functions.
Now if $f\in \Sigma(p)$, then $z/f \in \mathcal{H}$ and $(z/f)_{z=0}=1$. Therefore each function $z/f$ has the following Taylor expansion:
\beq\label{fp3eq3}
\frac{z}{f(z)}=1+b_{1}z+b_{2}z^2+\cdots,\quad z\in \ID.
\eeq
It is now natural to consider the following problems of maximizing the Yamashita functionals:
$$
\ds\max_ {f\in \Sigma(p)}\Delta(r,z/f)~ \mbox{ and} ~ \ds\max_ {f\in \mathcal{U}_{p}(\lambda)}\Delta(r,z/f).
$$
We answer the above problems in Section 3. Thereafter, we consider another problem that deals with finding the estimates of integral means
for the class $\Sigma(p)$ and its subclass  $\mathcal{U}_{p}(\lambda)$.
Now, consider $g\in \mathcal{H}$ and for such functions
define the integral means $L_1(r,g):=r^2 I_1(r,g)$ where
\beq\label{fp3eq2}
I_1(r,g)= \frac{1}{2\pi}\int_{-\pi}^{\pi}\frac{d\theta}{|g(re^{i\theta})|^2}, \quad 0<r\leq 1.
\eeq
We remark here that each $g\in \mathcal{H}$ has angular limits on the unit circle.
The above integral originated from a special case of integral means considered by Gromova and Vasilev in 2002 (see f.i. \cite{GV}).
The estimate for this integral has special applications in certain problems in fluid mechanics (compare \cite{AV-1, AV-2}). Recently Ponnusamy and Wirths
obtained sharp estimates of integral means for some subclasses of $\mathcal{A}$  (compare \cite{PW}) which settled one of the open problems of Gromova and Vasil’ev
described in \cite{GV}. We find sharp estimates for $L_1(r,f)$ whenever $f\in \Sigma(p)$ and its subclass  $\mathcal{U}_{p}(\lambda)$. These are also the contents
of Section ~3.

\section{Criteria for univalency and some preliminary results}
Let $F$ and $G$ be analytic in $\ID$. Now, a function $F$ is said to be
subordinate to $G$, written as $F\prec G$, if there exists a function $w$ analytic in $\ID$ with $w(0)=0$ and $|w(z)|<1$,
and such that $F(z)=G(w(z))$. If $G$ is univalent, then $F\prec G$ if and only if $F(0)=G(0)$ and $F(\ID)\subset G(\ID)$.
In the following theorem we prove a sufficient condition for univalence for $f\in\mathcal{A}(p)$.

\bthm
Let $f\in\mathcal{A}(p)$ with $f(z)/z \neq 0$ for $0<|z|<1$. If
$$
\left|\left(\frac{z}{f(z)}\right)''\right|\leq \left(\frac{1-p}{1+p}\right)^2,\quad z\in \ID,
$$
then $f$ is univalent in $\ID$.
\ethm
\bpf
Let $f\in\mathcal{A}(p)$ and since $f/z$ is nonvanishing in $\ID$, then $z/f$ is analytic in $\ID$ and has an expansion of the form (\ref{fp3eq3}).
From the given hypothesis,
$$
\left|-z^2\left(\frac{z}{f(z)}\right)''\right|\leq \left(\frac{1-p}{1+p}\right)^2|z|^2 <  \left(\frac{1-p}{1+p}\right)^2|z|,\quad z\in \ID.
$$
Therefore from the above inequality and applying the definition of subordination, we have
\beq\label{fp3eq15}
- z^2\left(\frac{z}{f(z)}\right)'' \prec z\left(\frac{1-p}{1+p}\right)^2.
\eeq
Now for $f\in\mathcal{A}(p)$, let
$$
p(z):=\frac{z}{f(z)}-z\left(\frac{z}{f(z)}\right)'=\left(\frac{z}{f(z)}\right)^2f'(z).
$$
Therefore, $p$ is analytic in $\ID$. Also it is a simple exercise to see that
$$p(z)=1+\sum_{n=1}^{\infty}(1-n)b_nz^n
$$
and
$$
p'(z)=-z \left(z/f(z)\right)'' \quad \mbox{i.e.} \quad zp'(z)= -z^2 \left(z/f(z)\right)''.
$$
Now, by (\ref{fp3eq15}) we get
$$
z p'(z) \prec z\left((1-p)/(1+p)\right)^2.
$$
By a consequence of a well known result of T. Suffridge (compare \cite[p. 76, Theorem 3.1d.]{MM}), we have
\beqq
&& p(z)\prec 1+ \frac{1}{2}\left(\frac{1-p}{1+p}\right)^2 z\\
&\mbox{i.e.}& \left(\frac{z}{f(z)}\right)^2f'(z) \prec 1+ \frac{1}{2}\left(\frac{1-p}{1+p}\right)^2 z\\
&\Rightarrow& \left|\left(\frac{z}{f(z)}\right)^2f'(z)-1\right|< \frac{1}{2}\left(\frac{1-p}{1+p}\right)^2 < \left(\frac{1-p}{1+p}\right)^2.
\eeqq
From the above inequality, we conclude that $f$ is univalent in $\ID$ by applying \cite[Theorem 1]{BF-1}.
\epf

We now move on to  present the absolute estimates
for the Dirichlet integrals $\Delta(r,f)$ and $\Delta(r,f/z)$ whenever $f\in\Sigma(p)$ and $0<r<p$.
We consider the sub-disc $\ID_p:=\{z:|z|<p\}\subsetneq \ID$.
We see that each $f\in \Sigma(p)$ has the Taylor expansion of the form
\beq\label{fp3eq8}
f(z)=z+\sum_{n=2}^{\infty}a_{n}z^{n},\quad z\in \ID_p.
\eeq

In 1962, Jenkins (\cite{jenk}) proved that if $f\in\Sigma(p)$ and has the form (\ref{fp3eq8}), then
\beq
\label{fp3eq9}
|a_n|\leq \frac{1+p^2+\cdots+p^{2n-2}}{p^{n-1}}=\frac{1-p^{2n}}{(1-p^2)p^{n-1}}, \quad n\geq 2.
\eeq
Equality holds in the above inequality for the function $k_p(z)=-pz/(z-p)(1-pz)$.
Now for $f\in\Sigma(p)$ we see that $f/z$ is analytic in $\ID_p$.
Therefore we consider the area problem for the functions $f/z$ whenever $f\in\Sigma(p)$ and $z\in\ID_p$. By using the Taylor expansion (\ref{fp3eq8}) for
$f\in \Sigma(p)$ and the Taylor coefficient estimate (\ref{fp3eq9}) for $|a_n|$, we have for $0<r<p$,
\beqq
\pi^{-1} \Delta(r,f/z)&=& \sum_{n=1}^{\infty}n|a_{n+1}|^2r^{2n}\\
                         &\leq& \sum_{n=1}^{\infty}n\left(\frac{1-p^{2n+2}}{(1-p^2)p^{n}}\right)^2r^{2n}\\
                         &=& \frac{p^2r^2}{(1-p^2)^2}\left(\frac{1}{(p^2-r^2)^2}-\frac{2}{(1-r^2)^2}+\frac{p^4}{(1-p^2r^2)^2}\right).
\eeqq
Equality holds in the above inequality for the function $k_p$.
Thus we infer that,
$$
\max_{f\in \Sigma(p)}\Delta(r,f/z)=\frac{\pi p^2r^2}{(1-p^2)^2}\left(\frac{1}{(p^2-r^2)^2}-\frac{2}{(1-r^2)^2}+\frac{p^4}{(1-p^2r^2)^2}\right)
$$
and the maximum is attained by the function $k_p$. Here we observe that
$$
\lim_{p\ra 1^{-}} \max_{f\in \Sigma(p)}\Delta(r,f/z)=2\pi r^2(r^2+2)(1-r^2)^{-4},
$$
which is same as the estimate obtained by Yamashita in (see \cite[p.435]{ Yama}) for $f\in \mathcal{S}$.
We also compute for $0<r<p$,
 \beqq
\max_{f\in \Sigma(p)}\Delta(r,f)= \Delta(r,k_p)=\frac{\pi p^2r^2}{(1-p^2)^2}\left(\frac{p^2}{(p^2-r^2)^2}-\frac{2}{(1-r^2)^2}
+\frac{p^2}{(1-p^2r^2)^2}\right).
\eeqq
As we pass through the limit $p\ra 1^{-}$, the right hand side of the above expression becomes $\pi r^2(r^4+4r^2+1)(1-r^2)^{-4}$.
We see that this estimate is same as the estimate obtained by Yamashita for the class $\mathcal{S}$ (Compare \cite[(4), p.436]{Yama}).

\section{Main Results}
We are now ready to state our first result after all the above discussion.

\bthm\label{fp3th1}
Let $f\in \Sigma(p)$ and $z/f$ have a Taylor expansion of the form $(\ref{fp3eq3})$ in $\ID$. Then for each $r\in (0,1]$, we have
$$
\max_{f\in \Sigma(p)}\Delta(r,z/f)=\pi r^2 \left((1/p+p)^2+2r^2\right)
$$
and the maximum is attained by the function $k_p(z)=-pz/(z-p)(1-pz)$.
\ethm

\bpf
Let $f\in \Sigma(p)$. We define $g(z)=cf(z)/(c+f(z))$ where $-c\notin f(\ID)$. It is easy to see that  $g\in \mathcal{S}$, and has the following Taylor expansion
$$
z/g(z)= z/c+z/f(z)=1+\left(b_{1}+1/c\right)z+b_{2}z^2+\cdots,\quad z\in \ID.
$$
Let $\ID^{*}\ni \xi=1/z$, then $F(\xi)=1/g(1/\xi)\in \Sigma$ and $F(\xi)$ takes the form
\[F(\xi)=\xi+(b_{1}+1/c)+ b_2/\xi+ b_3/\xi^2+\cdots,\quad \xi \in \ID^{*}.\]
Therefore, from the well known Gronwall's area theorem applied to the above function $F$, we have
\beq\label{fp3eq4}
\sum_{n=1}^{\infty}n|b_{n+1}|^2\leq1.
\eeq
We also observe from the expansion (\ref{fp3eq3}) and (\ref{fp3eq8}) that $b_1=-a_2$. Now using the inequality (\ref{fp3eq9}) for $n=2$, we have
 $|b_1|=|-a_2|\leq (1+p^2)/p$. Therefore by (\ref{fp3eq1}) and the expansion (\ref{fp3eq3}) we get
\beqq
\Delta(r,z/f)&=& \pi \sum_{n=1}^{\infty}n|b_n|^2r^{2n}\\
                        &=& \pi \left(|b_1|^2r^2+\sum_{n=1}^{\infty}(n+1)|b_{n+1}|^2r^{2(n+1)}\right)\\
                        &\leq& \pi \left(|b_1|^2r^2+\sum_{n=1}^{\infty}2n|b_{n+1}|^2r^{2(n+1)}\right)\\
                        &\leq&\pi \left((1/p+p)^2r^2+2 r^4\sum_{n=1}^{\infty}n|b_{n+1}|^2\right) \quad \left(\because |b_1|\leq(1+p^2)/p \right)\\
                        &\leq&\pi \left((1/p+p)^2r^2+2 r^4\right) ~~~~~~~~~~ \quad (\mbox{using}~ (\ref{fp3eq4}))\\
                        &=& \pi r^2\left((1/p+p)^2+2 r^2\right).
                        \eeqq
Equality holds in the above inequality for the function $k_p\in \Sigma(p)$. This can be easily seen if we observe that for this function $k_p$, we have
$b_1=-(1/p+ p),\quad b_2= 1$ and $b_n=0$ for $n\geq3$. Therefore, we conclude that
$$
\max_{f\in \Sigma(p)}\Delta(r,z/f)=\pi r^2 ((1/p+p)^2+2r^2).
$$
\epf

\br
As $p\ra 1^{-}$, the Dirichlet estimate in the above theorem is same as that of \cite[Theorem 1]{Yama}.
\er

In the following theorem we prove a sharp estimate for the integral mean $L_1(r,f)$ where $f\in \Sigma(p)$.
\bthm\label{fp3th2}
 Let $f\in \Sigma(p)$ and have the form $(\ref{fp3eq3})$. Then we have
 $$
 L_1(r,f)\leq 1+(1/p+p)^2r^2+r^4
 $$
 and the inequality is sharp.
 \ethm
 \bpf
  Let $f\in \Sigma(p)$ and have an expansion of the form (\ref{fp3eq3}). Then
  \beqq
  L_1(r,f)&:=& r^2I_1(r,f)\\
         &=& \frac{1}{2\pi}\int_{-\pi}^{\pi}\left|\frac{z}{f(z)}\right|^2d\theta= 1+\sum_{n=1}^{\infty}|b_n|^2r^{2n}\\
         &\leq& 1+|b_1|^2r^2+\sum_{n=2}^{\infty}(n-1)|b_n|^2r^{2n}\\
         &\leq& 1+(1/p+p)^2 r^2+r^4\sum_{n=1}^{\infty}n|b_{n+1}|^2r^{2n-2} \quad (\because |b_1|\leq(1+p^2)/p) \\
         &\leq& 1+(1/p+p)^2 r^2+r^4\sum_{n=1}^{\infty}n|b_{n+1}|^2 \quad (\mbox{since}~~  0<r\leq 1) \\
         &\leq& 1+(1/p+p)^2r^2+r^4 \quad (\mbox{by}~ (\ref{fp3eq4})).
 \eeqq
Equality holds in the above inequality for the function $k_p$.
 \epf

Now in a similar fashion, we can deduce a sharp estimate for the integral mean $L_1(r,f)$ where $f\in \mathcal{S}$. This is the
content of the following theorem.

\bthm\label{fp3th4}
Let $f\in \mathcal{S}$ and have the form $(\ref{fp3eq3})$. Then we have
$$
 L_1(r,f)\leq 1+4r^2+r^4
 $$
 and the result is sharp.
\ethm

\bpf
Let $f\in \mathcal{S}$ and have the expansion $(\ref{fp3eq3})$. We then have from \cite[Theorem 11, p.193. Vol.2]{goodm}
$$
\sum_{n=2}^{\infty}(n-1)|b_{n}|^2\leq 1.
$$
Here we note that $b_1=-a_2$ and from Bieberbach's theorem we know that $|a_2|\leq 2$,
with equality if and only if $f$ is a rotation of the Koebe function i.e. $f(z)=k_{\theta}(z)=z/(1-e^{i\theta}z)^2$ where $\theta$ is real.
Therefore from $(\ref{fp3eq2})$, we get
\beqq
  L_1(r,f)&:=& r^2I_1(r,f)\\
         &=& 1+\sum_{n=1}^{\infty}|b_n|^2r^{2n}\\
         &\leq& 1+|b_1|^2r^2+r^4\sum_{n=2}^{\infty}(n-1)|b_n|^2r^{2n-4}\\
         &\leq& 1+4 r^2+r^4\sum_{n=2}^{\infty}(n-1)|b_{n}|^2\quad \left(\because~ |b_1|=|-a_2|\leq 2~\mbox{and}~~ 0<r\leq1\right) \\
         &\leq& 1+4r^2+r^4 \quad \left(\because ~ \sum_{n=2}^{\infty}(n-1)|b_{n}|^2\leq1\right).
 \eeqq
Equality holds in the above inequality for the function $f=k_{\theta}$.
\epf

\br
Here we remark that as $p\ra1^{-}$, the integral mean in Theorem ~\ref{fp3th2} is
same as the integral mean that we obtain in Theorem~ \ref{fp3th4} for the class $\mathcal{S}$.
\er


We now move on to the class $\mathcal{U}_{p}(\lambda)$ and consider similar problems.
In doing so, we first prove the following Lemma which will be used to prove our main results for this function class.
We follow here the modified proof of \cite[Lemma 1]{OPW-2} and provide the details for the sake of completeness.
\blem\label{fp3lm1}
Let $f\in \mathcal{U}_{p}(\lambda)$ and have expansion of the form $(\ref{fp3eq3})$ for some $0<\lambda\leq 1$ and let $t\leq 2$. Then we have
$$
\sum_{n=2}^{\infty}n^t|b_n|^2r^{2n}\leq 2^t\lambda^2 \mu^2 r^4.
$$
\elem
\bpf
Suppose that $f\in \mathcal{U}_{p}(\lambda)$. Then we have (see \cite[Corollary 1]{BF-1})
$$
 |U_{f}(z)|\leq \lambda\mu |z|^{2},\quad z\in \ID,
 $$
 where
 $$
 U_f(z):=\left(z/f(z)\right)^2f'(z)-1=-z\left(z/f(z)\right)'+(z/f(z))-1.
 $$
 Now by the expansion (\ref{fp3eq3}) and the above inequality we get
 $$
 \left|\sum_{n=2}^{\infty}(n-1)b_nz^n\right|\leq \lambda\mu |z|^{2}.
 $$
 Therefore, for $z=re^{i\theta}$ and $0<r<1$,
 \beqq
 \sum_{n=2}^{\infty}(n-1)^2|b_n|^2r^{2n}&=& \frac{1}{2\pi}\int_{0}^{2\pi}\left|\sum_{n=2}^{\infty}(n-1)b_nz^n\right|^2\, d\theta \\
                                        &\leq& \lambda^2 \mu^2 r^4.
 \eeqq
 From the above inequality we get that for each $k\geq 2$, the inequality
 $$
 \sum_{n=2}^{k}(n-1)^2|b_n|^2r^{2n}\leq \lambda^2 \mu^2 r^4
 $$
 is true. Now, we consider these inequalities for $k=2,3,\cdots,N$, and multiply the $N$-th inequality by the factor $N^t/(N-1)^2,$
 and for $k=2,\cdots, N-1$, the $k$-th inequality by the factor
 $$
 \frac{k^t}{(k-1)^2}-\frac{(k+1)^t}{k^2}>0.
 $$
 Now after adding all these modified inequalities, we get in the left hand side of the inequality
 \beqq
 &&\sum_{k=2}^{N-1}\left(\left(\frac{k^t}{(k-1)^2}-\frac{(k+1)^t}{k^2}\right)\sum_{n=2}^{k}(n-1)^2|b_n|^2r^{2n}\right) +\frac{N^t}{(N-1)^2}
 \sum_{n=2}^{N}(n-1)^2|b_n|^2r^{2n} \\
 &=& \sum_{n=2}^{N-1}(n-1)^2|b_n|^2r^{2n}\frac{n^t}{(n-1)^2}+N^t|b_N|^2r^{2N} \\
 &=& \sum_{n=2}^{N}n^t|b_n|^2r^{2n}
 \eeqq
 and in the right hand side of the inequality, we get
 \beqq
 && \lambda^2 \mu^2 r^4 \left(\sum_{k=2}^{N-1}\left(\frac{k^t}{(k-1)^2}-\frac{(k+1)^t}{k^2}\right)+\frac{N^t}{(N-1)^2}\right) \\
 &=& 2^t\lambda^2 \mu^2 r^4.
 \eeqq
 As a result, we obtain the following inequality
 $$
 \sum_{n=2}^{N}n^t|b_n|^2r^{2n}\leq 2^t\lambda^2 \mu^2 r^4.
 $$
 Finally, letting $N\rightarrow \infty$, we have
 $$
 \sum_{n=2}^{\infty}n^t|b_n|^2r^{2n}\leq 2^t\lambda^2 \mu^2 r^4,
 $$
 which proves the lemma.
\epf

After plugging in $t=0$ in the above Lemma, we get
\beq\label{fp3eq5}
\sum_{n=2}^{\infty}|b_n|^2r^{2n}\leq \lambda^2 \mu^2 r^4
\eeq
and $t=1$, we get
\beq\label{fp3eq6}
\sum_{n=2}^{\infty}n|b_n|^2r^{2n}\leq 2\lambda^2 \mu^2 r^4.
\eeq

We are now in a position to state the following Theorem:
\bthm\label{fp3th3}
Let $f\in \mathcal{U}_{p}(\lambda)$ and have the form $(\ref{fp3eq3})$. Then we have
$$
\max_{f\in \mathcal{U}_{p}(\lambda)}\Delta(r,z/f)=\pi r^2\left((1/p+\lambda \mu p)^2+2\lambda^2 \mu^2 r^2\right)
$$
and
$$
L_1(r,f):=r^2 I_1(r,f)\leq 1+r^2(1/p+\lambda \mu p)^2+\lambda^2 \mu^2 r^4.
$$
The results are sharp for the function
\beq\label{fp3eq7}
f_p(z)=\frac{z}{1-\frac{z}{p}(1+\lambda\mu p^{2})+\lambda \mu z^{2}}, \quad z\in \ID.
\eeq
\ethm
\bpf
Let $f\in \mathcal{U}_{p}(\lambda)$.
Then we have, $|b_1|\leq 1/p+\lambda \mu p$ (Compare \cite[Theorem 5]{BF-1}). Now by (\ref{fp3eq1}) and (\ref{fp3eq6}) we have
\beqq
\Delta(r,z/f)&=& \pi \sum_{n=1}^{\infty}n|b_n|^2r^{2n} \\
             &=& \pi |b_1|^2r^2+\pi \sum_{n=2}^{\infty}n|b_n|^2r^{2n}\\
             &\leq& \pi r^2 (1/p+\lambda \mu p)^2+ \pi 2\lambda^2 \mu^2 r^4\\
             &=& \pi r^2\left((1/p+\lambda \mu p)^2+2\lambda^2 \mu^2 r^2\right).
\eeqq
To prove the sharpness assertion, we observe that $f_p\in \mathcal{U}_{p}(\lambda)$ and for the same
function $b_1=-(1/p+\lambda \mu p),\quad b_2= \lambda \mu$ and $b_n=0$ for $n\geq3$. Therefore it can be easily seen that equality occurs in the above inequality for the function $f_p$. Next we wish to prove the second part of the theorem. In order to do so, we compute using (\ref{fp3eq2}) and (\ref{fp3eq5}) that
\beqq
L_1(r,f):=r^2 I_1(r,f)&=&1+\sum_{n=1}^{\infty}|b_n|^2r^{2n}\\
                     &\leq & 1+r^2(1/p+\lambda \mu p)^2+\lambda^2 \mu^2 r^4.
\eeqq
We also see here that the above inequality is sharp for the function $f_p$.
\epf

Likewise for the analytic case, we also consider the classes of meromorphically convex ( abbreviated as {\it concave}) and meromorphically starlike
univalent functions in $\Sigma(p)$ which we denote by $Co(p)$ and $\Sigma^{*}(p,w_0)$ respectively. We clarify here that for $f\in Co(p)$, the set $\sphere \setminus f(\ID)$
is a compact convex set and for $f\in \Sigma^{*}(p,w_0)$, the compact set $\sphere \setminus f(\ID)$ is starlike with respect to a point $w_0\neq 0, \infty$.
The detailed discussion about these classes of functions can be found from \cite{APW1, APW2, AW, goodm}. Now we can deduce the following

\br
Let $f\in Co(p)\subsetneq \Sigma(p)$. Therefore,
$$
\max_{f\in Co(p)}\Delta(r,z/f)\leq\pi r^2 \left((1/p+p)^2+2r^2\right)
$$
and
$$
L_1(r,f)\leq 1+(1/p+p)^2r^2+r^4.
$$
As we know that $k_p\in Co(p)$, both the aforementioned results are sharp. Same conclusion can be drawn for $f\in \Sigma^*(p,w_0)$ as $k_p$ also belongs to the class $\Sigma^*(p,w_0)$ where $w_0 \in \left[-p/(1-p)^2, -p/(1+p)^2\right]$ (see \cite{BP-1}).
\er

\bigskip

{\bf Acknowledgement:} The authors thank Karl-Joachim Wirths for his
suggestions and careful reading of the manuscript.

\end{document}